\documentclass[a4paper,10pt, reqno]{amsart}
\usepackage{amsmath}
\usepackage{amssymb}
\usepackage{amsthm}
\usepackage{graphicx}
\usepackage{tabularx}
\usepackage{array}
\usepackage{physics}
\usepackage{cleveref}
\crefname{thm}{Theorem}{Theorems}
\crefname{lem}{Lemma}{Lemmas}
\crefname{prop}{Proposition}{Propositions}
\crefname{cor}{Corollary}{Corollaries}

\DeclareMathOperator{\sgn}{sgn}

\DeclareMathOperator{\CH}{conv}

\DeclareMathOperator{\CARD}{\#}

\newcommand{\R}{\mathbb{R}}

\newcommand{\Z}{\mathbb{Z}}
\newcommand{\PI}{\mathbb{Z}^{+}}

\newcommand{\origin}{o}

\newcommand{\SET}[2]{\left\{#1\;\middle|\; #2\right\}}

\newcommand{\SCB}[1][n]{\mathcal{K}^{#1}}

\newcommand{\floor}[1]{\left\lfloor #1 \right\rfloor}
\newtheorem{thm}{Theorem}

\newtheorem{lem}[thm]{Lemma}
\newtheorem{prop}[thm]{Proposition}

\theoremstyle{definition}

\newtheorem{rem}[thm]{Remark}

\begin{document}
\title{Covering functionals of convex polytopes with few vertices}
\author{Xia Li}
\email{lixia2016@nuc.edu.cn}
\author{Lingxu Meng}
\email{menglingxu@nuc.edu.cn}
\author{Senlin Wu}
\email{wusenlin@nuc.edu.cn}

\address{Department of Mathematics, North University of China, 030051 Taiyuan,
  China}



\begin{abstract}
  By using elementary yet interesting observations and refining techniques used
  in a recent work by Fei Xue et al., we present new upper bounds for covering
  functionals of convex polytopes in $\R^n$ with few vertices. In
  these estimations, no information other than the number of vertices of the
  convex polytope is used.
\end{abstract}

\maketitle

\section{Introduction}
\label{sec:introduction}

Hadwiger's covering conjecture is a long-standing open problem from convex and
discrete geometry which asserts that each \emph{convex body} (i.e., each compact
convex set having interior points) in $\R^n$ can be covered by at most $2^n$ of
its smaller homothetic copies. Despite all the efforts made, this conjecture is
still open even when $n=3$. See, e.g.,
\cite{Boltyanski-Martini-Soltan1997,Martini-Soltan1999,Brass-Moser-Pach2005,Bezdek2010,Bezdek-Khan2018}
and references therein for the history, classical results, as well as recent
progress about this conjecture. To attack Hadwiger's covering conjecture,
Chuanming Zong introduced the first computer-based quantitative program (cf.
\cite{Zong2010}), which is promising if the conjecture admits an affirmative
answer.

One major part of Zong's program is to get an efficient estimation of the upper
bound of covering functionals of each member of a properly chosen
$\varepsilon$-net of the space $\SCB$ of convex bodies in $\R^n$ endowed with
the Banach-Mazur metric, where for each compact convex set $K$ and each
$m\in\PI$, the \emph{covering functional $\Gamma_m(K)$ of $K$ with respect to
  $m$} is defined by
\begin{displaymath}
  \Gamma_m(K)=\inf\SET{\gamma>0}{\exists C\subset\R^n\text{ with }\CARD C=m\text{ s.t.
    }K\subseteq C+\gamma K}.
\end{displaymath}
Since each convex body can be approximated by convex polytopes, such an
$\varepsilon$-net can be chosen from the set of convex polytopes in $\SCB$.

In a recent work \cite{Wu-He2019}, Senlin Wu and Chan He provided a method to
obtain upper bounds for covering functionals of convex polytopes. In particular,
it is shown that when $K\in \SCB$ is a convex polytope with vertex set $V$ and
$\CARD V\leq 2^n$, then by \cite[Corollary 4]{Wu-He2019},
\begin{displaymath}
  \Gamma_{2^n}(K)\leq 1-\lambda(K),
\end{displaymath}
where $\lambda(K)=\inf\SET{\lambda(x,K)}{x\in K}$ and 
\begin{displaymath}
  \lambda(x,K)=\sup\SET{\lambda\in [0,1]}{\exists v\in V\text{ and }y\in K\text{
      s.t. }x=\lambda v+(1-\lambda)y},~\forall x\in K.
\end{displaymath}
The functional $\lambda(K)$ satisfies (cf. \cite[Lemma 10
and Theorem 12]{Wu-He2019})
\begin{displaymath}
  \frac{1}{n+1}\leq\lambda(K)\leq\frac{1}{2},
\end{displaymath}
and is not easy to determine in general. When there is no information on the
value of $\lambda(K)$, we only have the following general estimation:
\begin{equation}
  \label{eq:general-estimation}
  \Gamma_{2^n}(K)\leq 1-\frac{1}{n+1}.
\end{equation}

\bigskip

In this paper, by applying and refining techniques in \cite{Xue-Lian-Zhang2021+}
by Fei Xue et al., we provide new upper bounds for $\Gamma_{2^n}(K)$ when $K$ is
a convex polytope in $\SCB$ having at most $2^n$ vertices without using any
knowledge about $\lambda(K)$. For this purpose, we collect several technical
lemmas in Section \ref{sec:auxiliary-lemmas}. In Section \ref{sec:Kp}, we
estimate covering functionals of the $n$-dimensional $\ell_p~(p\geq 1)$ ball
$K_n^p$ and the portion of $K_n^{p}$ in the nonnegative orthant, which is
denoted by $K_n^{p*}$. More precisely,
\begin{gather*}
  K_n^p=\SET{(x_1,\dots,x_n)\in\R^n}{\sum\limits_{i\in[n]}|x|_i^p\leq 1},\\
  K_n^{p*}=\SET{(x_1,\dots,x_n)\in\R^n}{\sum\limits_{i\in[n]}x_i^p\leq 1,~x_i\geq 0,\forall i\in [n]}.
\end{gather*}
Here, we used the shorthand notation $[m]:=\SET{i\in\PI}{1\leq i\leq m}$ for
each $m\in\PI$. In Section \ref{sec:few-vertices}, we provide upper bounds for
$\Gamma_{2^n}(K)$ in each of the following two cases: $K$ is a convex polytope
with $m~(m\geq n+1)$ vertices; $K$ is a centrally symmetric convex polytope
having $2m~(m\geq n)$ vertices.

Throughout this paper, the dimension $n$ of the underlying space is at least
$3$.

\section{Auxiliary Lemmas}
\label{sec:auxiliary-lemmas}

In this section, $f$ and $g$ denotes the functions
\begin{align*}
  f:~(0,\infty)&\to \R &g:~(0,\infty)&\to\R\\
  x&\mapsto \frac{(1+x)^{1+x}}{x^x},&x&\mapsto \frac{2^x(1+x)^{1+x}}{x^x}, 
\end{align*}
respectively. One can easily verify that both $f$ and $g$ are strictly
increasing on $(0,+\infty)$, and that $\lim\limits_{x\to
  0+}f(x)=\lim\limits_{x\to 0+}g(x)=1$. Note also that $f(1)=4$ and $g(1)=8$.

For each $t\in (1,\infty)$, let $a(t)$ and $b(t)$ be the solution to the
equation $f(x)=t$ and $g(x)=t$, respectively. Since $f$ and $g$ are both
strictly increasing, $a(t)$ and $b(t)$ are also strictly increasing. Numerical
calculation shows that $a(2)\approx 0.293815$ and $b(2)\approx 0.205597$.

The main result in \cite{Robbins1955} shows that 
\begin{displaymath}
  n!=\sqrt{2\pi}n^{n+\frac{1}{2}}e^{-n}e^{r_n},
\end{displaymath}
where
\begin{displaymath}
  \frac{1}{12n+1}<r_n<\frac{1}{12n}.
\end{displaymath}
It follows that, for each $k\in\PI$,
\begin{align}
  \binom{n+k}{n}=\frac{(n+k)!}{n!k!}&=\frac{\sqrt{2\pi}(n+k)^{n+k+\frac{1}{2}}e^{-(n+k)}e^{r_{n+k}}}{\sqrt{2\pi}n^{n+\frac{1}{2}}e^{-n}e^{r_n}\sqrt{2\pi}k^{k+\frac{1}{2}}e^{-k}e^{r_k}}\nonumber\\
                                    &=\qty(\frac{n+k}{2\pi nk})^{\frac{1}{2}}\qty[\frac{\qty(1+\frac{k}{n})^{1+\frac{k}{n}}}{\qty(\frac{k}{n})^{\frac{k}{n}}}]^ne^{r_{n+k}-r_n-r_k}.\label{eq:first-binom}
\end{align}
Since
\begin{displaymath}
  r_{n+k}-r_n-r_k<\frac{1}{12(n+k)}-\frac{1}{12n+1}-\frac{1}{12k+1}<0,
\end{displaymath}
we have
\begin{equation}
  \label{eq:estimation-binom}
  \binom{n+k}{n}<\qty(\frac{n+k}{2\pi nk})^{\frac{1}{2}}\qty[\frac{\qty(1+\frac{k}{n})^{1+\frac{k}{n}}}{\qty(\frac{k}{n})^{\frac{k}{n}}}]^n=\qty(\frac{n+k}{2\pi nk})^{\frac{1}{2}}\qty[f\qty(\frac{k}{n})]^n.
\end{equation}

For each $t\in(1,2]$ and each $n\in \PI$, let $k(n,t)$ and $l(n,t)$ be the nonnegative integers satisfying
\begin{displaymath}
  \binom{n+k(n,t)}{n}\leq t^n<\binom{n+k(n,t)+1}{n}
\end{displaymath}
and
\begin{displaymath}
  2^{l(n,t)}\binom{n+l(n,t)}{n}\leq t^n<2^{l(n,t)+1}\binom{n+l(n,t)+1}{n},
\end{displaymath}
respectively.

\begin{lem}
  \label{lem:estimations-k(n)}
  Let $t\in(1,2]$ and $n\in \PI$. We have
  \begin{enumerate}
  \item[(a)] $k(n,t)\geq \floor{a(t)n}$;
  \item[(b)] $k(n,t)>0$ if and only if $n+1\leq t^n$;
  \item[(c)] $l(n,t)\geq \floor{b(t)n}$;
  \item[(d)] $l(n,t)>0$ if and only if $2(n+1)\leq t^n$.
  \end{enumerate}  
\end{lem}
\begin{proof}
  (a). The case when $\floor{a(t)n}=0$ is clear. When $\floor{a(t)n}\geq 1$, by
  \cref{eq:estimation-binom}, we have
  \begin{align*}
    \binom{n+\floor{a(t)n}}{n}&< \qty(\frac{n+\floor{a(t)n}}{2\pi n\floor{a(t)n}})^{\frac{1}{2}}\qty[f\qty(\frac{\floor{a(t)n}}{n})]^n\\
                              &\leq \qty(\frac{n+a(t)n}{2\pi n})^{\frac{1}{2}}\qty[f\qty(\frac{a(t)n}{n})]^n\\
                              &=\qty(\frac{1+a(t)}{2\pi})^{\frac{1}{2}}t^n\\
    &\leq \qty(\frac{1+a(2)}{2\pi})^{\frac{1}{2}}t^n<t^n.
  \end{align*}
  Therefore
  \begin{displaymath}
    \binom{n+\floor{a(t)n}}{n}<\binom{n+k(n,t)+1}{n},
  \end{displaymath}
  which shows that $\floor{a(t)n}<k(n,t)+1$. Thus $k(n,t)\geq \floor{a(t)n}$.

  (b). It sufficies to observe that
  \begin{displaymath}
    k(n,t)=0\Longleftrightarrow 1=\binom{n}{n}\leq t^n<\binom{n+1}{n}=n+1.
  \end{displaymath}

  (c). We only need to consider the case when $\floor{b(t)n}\geq 1$. By
  \cref{eq:estimation-binom}, we have
  \begin{align*}
    2^{\floor{b(t)n}}\binom{n+\floor{b(t)n}}{n}&<
                                                 \qty(\frac{n+\floor{b(t)n}}{2\pi
                                                 n\floor{b(t)n}})^{\frac{1}{2}}\qty[g\qty(\frac{\floor{b(t)n}}{n})]^n\\
                                               &\leq\qty(\frac{n+b(t)n}{2\pi
                                                 n})^{\frac{1}{2}}[g\qty(b(t))]^n\\
                                               &\leq \qty(\frac{1+b(t)}{2\pi})^{\frac{1}{2}}t^n<t^n.
  \end{align*}
  It follows that
  \begin{displaymath}
    2^{\floor{b(t)n}}\binom{n+\floor{b(t))n}}{n}<2^{l(n,t)+1}\binom{n+l(n,t)+1}{n}.
  \end{displaymath}
  Hence $l(n,t)\geq \floor{b(t)n}$.

  (d). As for (b), we only need to observe that
  \begin{displaymath}
    l(n,t)=0\Longleftrightarrow 1=\binom{n}{n}\leq t^n<2\binom{n+1}{n}=2(n+1).\qedhere
  \end{displaymath}
\end{proof}

\begin{lem}
  \label{lem:lp}
  If either $x\geq a\geq 0$ or $x\leq a\leq 0$, then $|x-a|^p\leq |x|^p-|a|^p$
  holds for each $p\in[1,\infty)$.
\end{lem}
\begin{proof}
  First, suppose that $a\geq 0$. Set
  \begin{displaymath}
    h(x)=x^p-(x-a)^p,~\forall x\in[a,\infty).
  \end{displaymath}
  Then $h'(x)=px^{p-1}-p(x-a)^{p-1}\geq 0$. Thus $h$ is increasing on
  $[a,\infty)$. It follows that $h(x)\geq h(a)=a^p,~\forall x\in[a,\infty)$.

  When $x\leq a\leq 0$, we  have
  \begin{displaymath}
    |x-a|^p=((-x)-(-a))^p\leq (-x)^p-(-a)^p=|x|^p-|a|^p.\qedhere
  \end{displaymath}
\end{proof}

For $n,k\in \PI$, we put
\begin{gather*}
  M_1(n,k)=\SET{(x_1,\dots,x_n)\in\Z^n}{\sum\limits_{i\in[n]} x_i\leq k,~x_i\geq
    0,~\forall i\in[n]},\\
  M_2(n,k)=\SET{(x_1,\dots,x_n)\in\Z^n}{\sum\limits_{i\in[n]} |x_i|\leq k}.
\end{gather*}
It can be verified that (see e.g., \cite[Proposition 2.1]{Xue-Lian-Zhang2021+})
\begin{equation}
  \label{eq:card-M1}
  \CARD M_1(n,k)=\binom{n+k}{n}.
\end{equation}
And it is known that (see e.g., \cite[p. 182]{Polya-Szego1998} or \cite{Betke-Henk1993})
\begin{displaymath}
  \CARD M_2(n,k)=\sum\limits_{i=0}^n2^{n-i}\binom{n}{i}\binom{k}{n-i},
\end{displaymath}
which implies that (cf. also \cite[p. 7]{Xue-Lian-Zhang2021+})
\begin{align}
  \CARD M_2(n,k)=\sum\limits_{i=n-k}^{n}2^{n-i}\binom{n}{i}\binom{k}{n-i}&\leq 2^k\sum\limits_{i=n-k}^n\binom{n}{i}\binom{k}{n-i}\nonumber\\
                                                                   &=2^k\sum\limits_{i=0}^n\binom{n}{i}\binom{k}{n-i}=2^k\binom{n+k}{n}.\label{eq:card-M2}
\end{align}

\section{Covering functionals of $K_n^{p*}$ and $K_n^p$}
\label{sec:Kp}

\cref{lem:inclusion-Knp-star} and \cref{lem:p-ball} are appeared in
\cite{Xue-Lian-Zhang2021+}. For readers' convenience we provide new proofs
without using induction.

For each $n\in\PI$ and $p\in [1,\infty)$, put
\begin{displaymath}
  \overline{K_n^{p*}}=\SET{(x_1,x_2,\dots,x_n)\in\R^n}{\sum\limits_{i\in[n]}x_i^p\leq
    n,~x_i\geq 0,~\forall i\in[n]}.
\end{displaymath}

\begin{lem}
  \label{lem:inclusion-Knp-star}
  For each pair $n,k\in \PI$ and each $p\geq 1$, we have
  \begin{displaymath}
    \qty(\frac{n+k}{n})^{\frac{1}{p}}\overline{K_n^{p*}}\subseteq \overline{K_n^{p*}}+M_1(n,k).
  \end{displaymath}  
\end{lem}
\begin{proof}
  Let $(z_1,\dots,z_n)$ be an arbitrary point in
  $\qty(\frac{n+k}{n})^{\frac{1}{p}}\overline{K_n^{p*}}$. Then
  \begin{displaymath}
    z_i\geq 0,~\forall i\in[n]\qqtext{and} \sum\limits_{i\in[n]} z_i^p\leq n+k.
  \end{displaymath}
  Set $n_0:=\CARD\SET{z_i}{z_i\geq 1,~i\in[n]}$. Without loss of generality,
  assume that $z_1,\dots,z_{n_0}\geq 1$. We distinguish three cases.

  {\bfseries Case 1:} $n_0\geq k$. In this situation, we have
  \begin{displaymath}
    (z_1,\dots,z_n)=(z_1-1,\dots,z_k-1,z_{k+1},\dots,z_n)+(\underbrace{1,\dots,1}_{k},0,\dots,0).
  \end{displaymath}
  From \cref{lem:lp} it follows that
  \begin{displaymath}
    \sum\limits_{i\in[k]}(z_i-1)^p+\sum\limits_{i\in[n]\setminus[k]}z_i^p\leq
    \sum\limits_{i\in[k]}(z_i^p-1)+\sum\limits_{i\in[n]\setminus[k]}z_i^p\leq (n+k)-k=n.
  \end{displaymath}
  Thus
  \begin{equation}
    \label{eq:K-n-p}
    (z_1,\dots,z_n)\in \overline{K_n^{p*}}+M_1(n,k).
  \end{equation}

  {\bfseries Case 2: }$n_0<k$, $\sum\limits_{i\in[n_0]}\floor{z_i}\geq k$. There exist integers $a_1,\dots,a_{n_0}\geq 0$ such that
  \begin{displaymath}
    a_i\leq z_i,~\forall i\in[n_0]\qqtext{and}\sum\limits_{i\in[n_0]}a_i=k.
  \end{displaymath}
  Thus
  \begin{displaymath}
    (z_1,\dots,z_n)=(z_1-a_1,\dots,z_{n_0}-a_{n_0},z_{n_0+1},\dots,z_n)+(a_1,\dots,a_{n_0},0,\dots,0).
  \end{displaymath}
  By \cref{lem:lp}, 
  \begin{displaymath}
    \sum\limits_{\in[n_0]}(z_i-a_i)^p+\sum\limits_{i\in[n]\setminus[n_0]}z_i^p\leq
    \sum\limits_{i\in[n_0]}(z_i^p-a_i^p)+\sum\limits_{i\in[n]\setminus[n_0]}z_i^p\leq (n+k)-k=n.
  \end{displaymath}
  Again, we have \cref{eq:K-n-p}.

  {\bfseries Case 3: }$n_0<k$ and $\sum\limits_{i\in[n_0]}\floor{z_i}<k$. We have the following decomposition:
  \begin{displaymath}
    (z_1,\dots,z_n)=(z_1-\floor{z_1},\dots,z_{n_0}-\floor{z_{n_0}},z_{n_0+1},\dots,z_n)+(\floor{z_1},\dots,\floor{z_{n_0}},0,\dots,0).
  \end{displaymath}
  Since $z_i-\floor{z_i}<1,~\forall i\in[n_0]$ and $z_i<1,~\forall i\in[n]\setminus[n_0]$, we have
  \begin{displaymath}
    \sum\limits_{i\in[n_0]}(z_i-\floor{z_i})^p+\sum\limits_{i\in[n]\setminus[n_0]}z_i^p\leq n,
  \end{displaymath}
  which implies \cref{eq:K-n-p}.
\end{proof}

For each $p\geq 1$ and each $n\in\PI$, let
\begin{displaymath}
  \overline{K_n^p}=\SET{(x_1,x_2,\dots,x_n)\in\R^n}{\sum\limits_{i=1}^n|x_i|^p\leq  n,~x_i\geq 0,~\forall i=1,\dots,n}.
\end{displaymath}

\begin{lem}
  \label{lem:p-ball}
  For each pair $n,k\in\PI$ and each $p\geq 1$, we have
  \begin{displaymath}
    \qty(\frac{n+k}{n})^{\frac{1}{p}}\overline{K_n^p}\subseteq \overline{K_n^p}+M_2(n,k).
  \end{displaymath}
\end{lem}
\begin{proof}
  Let $(z_1,\dots,z_n)$ be an arbitrary point in
  $\qty(\frac{n+k}{n})^{\frac{1}{p}}\overline{K_n^p}$. Then
  \begin{displaymath}
    \sum\limits_{i\in[n]}|z_i|^p\leq n+k. 
  \end{displaymath}
  Set $n_0=\CARD\SET{z_i}{|z_i|\geq 1,~i\in[n]}$. Without loss of generality,
  assume that $|z_1|,\dots,|z_{n_0}|\geq 1$. We distinguish three cases.

  {\bfseries Case 1: }$n_0\geq k$. Clearly,
  \begin{align*}
    (z_1,\dots,z_n)=&(z_1-\sgn z_1,\dots,z_{k}-\sgn{z_{k}},z_{k+1},\dots,z_n)\\
                    &+(\sgn z_1,\dots,\sgn{z_{k}},0,\dots,0).
  \end{align*}
  By \cref{lem:lp}, we have
  \begin{displaymath}
    \sum\limits_{i\in[k]}|z_i-\sgn{z_i}|^p+\sum\limits_{i\in[n]\setminus[k]}|z_i|^p\leq
    \sum\limits_{i\in[k]}(|z_i|^p-|\sgn z_i|^p)+\sum\limits_{i\in[n]\setminus[k]}|z_i|^p\leq (n+k)-k=n.
  \end{displaymath}
  Thus
  \begin{equation}
    \label{eq:lp-ball}
    (z_1,\dots,z_n)\in \overline{K_n^p}+M_2(n,k).
  \end{equation}

  {\bfseries Case 2: }$n_0<k$ and $\sum\limits_{i\in[n_0]}\floor{|z_i|}\geq k$.
  There exist nonnegative integers $a_1,\dots,a_{n_0}$ such that
  \begin{displaymath}
    a_i\leq |z_i|,~\forall i\in[n_0]\qqtext{and}\sum\limits_{i\in[n_0]}a_i=k.
  \end{displaymath}
  Then
  \begin{align*}
    (z_1,\dots,z_n)=&(z_1-a_1\sgn{z_1},\dots,z_{n_0}-a_{n_0}\sgn{z_{n_0}},z_{n_0+1},\dots,z_n)\\
                    &+(a_1\sgn{z_1},\dots,a_{n_0}\sgn{z_{n_0}},0,\dots,0).
  \end{align*}
  Applying \cref{lem:lp} again, we obtain
  \begin{displaymath}
    \sum\limits_{i\in[n_0]}|z_i-a_i\sgn{z_i}|^p+\sum\limits_{i\in[n]\setminus[n_0]}|z_i|^p\leq
    \sum\limits_{i\in[n_0]}(|z_i|^p-|a_i\sgn z_i|^p)+\sum\limits_{i\in[n]\setminus[n_0]}|z_i|^p\leq n.
  \end{displaymath}
  Hence we have \cref{eq:lp-ball} again.

  {\bfseries Case 3: }$n_0<k$ and $\sum\limits_{i\in[n_0]}\floor{|z_i|}<k$. We have the decomposition:
  \begin{align*}
    (z_1,\dots,z_n)=&(z_1-\floor{|z_1|}\sgn{z_1},\dots,z_{n_0}-\floor{|z_{n_0}|}\sgn{z_{n_0}},z_{n_0+1},\dots,z_n)\\
                    &+(\floor{|z_1|}\sgn{z_1},\dots,\floor{|z_{n_0}|}\sgn{z_{n_0}},0,\dots,0).
  \end{align*}
  Moreover,
  \begin{displaymath}
    \sum\limits_{i\in[n_0]}|z_1-\floor{|z_1|}\sgn{z_1}|^p+\sum\limits_{i\in[n]\setminus[n_0]}|z_i|^p\leq
    n\qqtext{and}\sum\limits_{i\in[n_0]}|\floor{|z_i|}\sgn{z_i}|\leq k.
  \end{displaymath}
  It follows that \cref{eq:lp-ball} holds.
\end{proof}

\begin{prop}
  \label{prop:Kpn}
  Let $p\geq 1$ and $t\in(1,2]$. We have
  \begin{gather*}
    \Gamma_{t^n}(K_n^{p*})\leq \qty(\frac{n}{n+k(n,t)})^{\frac{1}{p}}\leq
    \qty(\frac{n}{n+\floor{a(t)n}})^{\frac{1}{p}},\\
    \Gamma_{t^n}(K_n^{p})\leq\qty(\frac{n}{n+l(n,t)})^{\frac{1}{p}}\leq
    \qty(\frac{n}{n+\floor{b(t)n}})^{\frac{1}{p}}.
  \end{gather*}
\end{prop}
\begin{proof}
  By \cref{lem:inclusion-Knp-star} we have
  \begin{displaymath}
    \Gamma_{t^n}(K_n^{p*})\leq \Gamma_{\binom{n+k(n,t)}{n}}(K_n^{p*})=\Gamma_{\binom{n+k(n,t)}{n}}(\overline{K_n^{p*}})\leq \qty(\frac{n}{n+k(n,t)})^{\frac{1}{p}}
  \end{displaymath}
  Applying (a) in \cref{lem:estimations-k(n)}, we obtain
  \begin{displaymath}
    \Gamma_{t^n}(K_n^{p*})\leq \qty(\frac{n}{n+\floor{a(t)n}})^{\frac{1}{p}}.
  \end{displaymath}
  It can be easily verified that
  \begin{displaymath}
    \Gamma_{t^n}(K_n^{p})\leq \Gamma_{2^{l(n,t)}\binom{n+l(n,t)}{n}}(K_n^{p})\leq\Gamma_{\CARD M_2(n,l(n,t))}(K_n^{p})=\Gamma_{\CARD M_2(n,l(n,t))}(\overline{K_n^{p}}).
  \end{displaymath}
  By \cref{lem:p-ball} and (c) in \cref{lem:estimations-k(n)}, we have
  \begin{displaymath}
    \Gamma_{t^n}(K_n^{p})\leq \qty(\frac{n}{n+l(n,t)})^{\frac{1}{p}} \leq \qty(\frac{n}{n+\floor{b(t)n}})^{\frac{1}{p}}.\qedhere
  \end{displaymath}
\end{proof}

\begin{rem}
  For $n\geq 1$, $n+1\leq 2^n$. So $k(n,2)>0$ by (b) in \cref{lem:estimations-k(n)}. Then, 
  \begin{displaymath}
    \Gamma_{2^n}(K_n^{p*})\leq \qty(\frac{n}{n+k(n,2)})^{\frac{1}{p}}\leq \qty(\frac{n}{n+1})^{\frac{1}{p}}<1.
  \end{displaymath}
  Similarly, if $n\geq 3$ (or equivalently, if $n+1\leq 2^{n-1}$), then
  $l(n,2)>0$ by (d) in \cref{lem:estimations-k(n)}. In this case we have
  \begin{displaymath}
    \Gamma_{2^n}(K_n^p)\leq \qty(\frac{n}{n+l(n,2)})^{\frac{1}{p}}\leq \qty(\frac{n}{n+1})^{\frac{1}{p}}<1.
  \end{displaymath}
\end{rem}

Let $K_1$ and $K_2$ be two convex bodies in $\R^n$. The multiplicative
\emph{Banach-Mazur distance} $d_{BM}^M(K_1,K_2)$ is defined by
\begin{displaymath}
  d_{BM}^M(K_{1},K_{2}):=\min\limits_{x\in\R^n, T\in\mathcal{A}^{n}} \SET{\gamma\geq 1}{K_{1}\subseteq
    T(K_{2})\subseteq\gamma K_{1}+x},
\end{displaymath}
where $\mathcal{A}^n$ is the set of non-singular affine transformations from
$\R^n$ to $\R^n$. When $K_1$ and $K_2$ are both centrally symmetric with respect
to the origin, $d_{BM}^M(K_1,K_2)$ is equal to the Banach-Mazur distance between
the Banach spaces having $K_1$ and $K_2$ as unit balls. Moreover,
\cite[Proposition 37.6]{Jaegermann1989} shows that
\begin{displaymath}
  d_{BM}^M(K_n^p,K_n^\infty)=n^{\frac{1}{p}},~\forall p\geq 2.
\end{displaymath}
This, together with the following fact proved in the proof of Theorem A in \cite{Zong2010}
\begin{displaymath}
  |\Gamma_{m}(K_1)-\Gamma_{m}(K_2)|\leq d_{BM}^M(K_1,K_2)-1,
\end{displaymath}
implies that 
\begin{displaymath}
  \Gamma_{2^n}(K_n^{p})\leq\min\qty{\qty(\frac{n}{n+\floor{b(2)n}})^{\frac{1}{p}},n^{\frac{1}{p}}-\frac{1}{2}},~\forall
  p\geq 2.
\end{displaymath}
For $1\leq p<2$, $n^{\frac{1}{p}}-\frac{1}{2}\geq 1$ since $n\geq 3$. Thus the
inequality above is also valid for the case when $p\in[1,2)$. Observe that
$\qty(\frac{n}{n+\floor{b(2)n}})^{\frac{1}{p}}$ is increasing and
$n^{\frac{1}{p}}-\frac{1}{2}$ is strictly decreasing with respect to $p$. Let
$p(n)$ be the unique solution to
\begin{displaymath}
  \qty(\frac{n}{n+\floor{b(2)n}})^{\frac{1}{p}}=n^{\frac{1}{p}}-\frac{1}{2}.
\end{displaymath}
We have
\begin{displaymath}
  \Gamma_{2^n}(K_n^{p})\leq\qty(\frac{n}{n+\floor{b(2)n}})^{\frac{1}{p(n)}},~\forall
  p\geq 1,
\end{displaymath}
which provides a uniform upper bound for $\Gamma_{2^n}(K_n^p)$ when $n$ is fixed
(see \cite[Theorem 2]{Zong2010} for a result of this type when $n=3$). It can
be shown that
\begin{displaymath}
  \frac{\ln n}{\ln\qty(\frac{3}{2})}\leq p(n)\leq \frac{\ln n}{\ln\qty(\frac{1}{2}+\frac{1}{1+b(2)})}.
\end{displaymath}

\section{Covering functionals of convex polytopes with few vertices}
\label{sec:few-vertices}
\begin{lem}\label{lem:linear-operator}
  Let $K\subseteq \R^m$ be a compact convex set, $T:~\R^m\to\R^n$ be an affine
  transformation. Then, for each $l\in\PI$, we have
  \begin{displaymath}
    \Gamma_l(T(K))\leq \Gamma_l(K).
  \end{displaymath}
\end{lem}
\begin{proof}
  Since $\Gamma_l(K)$ is translation invariant, we only need to consider the
  case when $T$ is linear. By the definition of $\Gamma_l(K)$, there exists a
  set $C\subseteq\R^m$ with $\CARD C=l$ such that $K\subseteq C+\Gamma_l(K)K$.
  Since $T$ is linear, we have
  \begin{displaymath}
    T(K)\subseteq T\qty(C+\Gamma_l(K)K)=T(C)+\Gamma_{l}(K)T(K).
  \end{displaymath}
  Since $\CARD T(C)\leq \CARD C=l$, we have
  \begin{displaymath}
    \Gamma_l(T(K))\leq \Gamma_{\CARD{T(C)}}(T(K))\leq \Gamma_l(K).\qedhere
  \end{displaymath}
\end{proof}

For each $m\in \PI$, we denote by $\Delta_m$ an $m$-dimensional simplex. I.e.,
$\Delta_m$ is the convex hull of $m+1$ affinely independent points. One can
easily verify that $\Delta_m$ is affinely equivalent to $K_m^{1*}$. A compact
convex set $K$ is said to be an \emph{$m$-dimensional crosspolytope} if there
exist $m$ linearly independent vectors $v_1,\dots,v_m$ such that
\begin{displaymath}
  K=\CH\qty{\pm v_1,\dots,\pm v_m}.
\end{displaymath}
Clearly, any $m$-dimensional crosspolytope is affinely equivalent to $K_m^1$.

\begin{thm}
  \label{thm:simplex-crosspolytope}
  Let $K\subseteq\R^n$ be a convex polytope.
  \begin{enumerate}
  \item[(a).] If $K$ has $m+1$ vertices, then $\Gamma_l(K)\leq \Gamma_l(K_m^{1*}),~\forall l\in\PI$.
  \item[(b).] If $K$ is centrally symmetric and has $2m$ vertices, then
    $\Gamma_l(K)\leq \Gamma_l(K_m^1)$ holds for any $l\in \PI$.
  \end{enumerate}  
\end{thm}
\begin{proof}
  (a). Let $v_0,v_1,\dots,v_m$ be the vertices of $K$. Denote by $\origin_m$ the
  origin of $\R^m$. For each $i\in[m]$, we denote by $e_i^m$ the $i$-th
  member of the standard basis of $\R^m$. We identify each ordered pair
  $(a,b)\in \R^n\times \R^m$, where $a=(\alpha_1,\dots,\alpha_n)$ and
  $b=(\beta_1,\dots,\beta_m)$, with the point
  $(\alpha_1,\dots,\alpha_n,\beta_1,\dots,\beta_m)\in\R^{n+m}$. Put
  \begin{displaymath}
    u_0=(v_0,\origin_m)\qqtext{and} u_i=(v_i,e_i^m),~\forall i\in[m].
  \end{displaymath}
  We claim that $u_0, u_1,\dots,u_m$ are affinely independent. Indeed, if
  $\lambda_0,\lambda_1,\dots,\lambda_m\in\R$ statisfy that
  \begin{displaymath}
    \sum\limits_{i\in[m]}\lambda_i=0\qqtext{and} \sum\limits_{i\in[m]}\lambda_iu_i=\origin_{n+m},
  \end{displaymath}
  then $\sum\limits_{i\in[m]} \lambda_ie_i^m=\origin_m$, which shows that
  $\lambda_i=0,~\forall i\in[m]$. Hence $\lambda_i=0$ holds for each $i\in[m]\cup\qty{0}$. Therefore, $L=\CH\qty{u_0,u_1,\dots,u_m}$ is an
  $m$-dimensional simplex. Since $\Gamma_l(\cdot)$ is an affine invariant,
  $\Gamma_l(L)=\Gamma_l(K_m^{1*})$. Let
  \begin{equation}
    \label{eq:map}
    T:~\R^{m+n}\to \R^n,\quad (\alpha_1,\dots,\alpha_{m+n})\mapsto (\alpha_1,\dots,\alpha_n).
  \end{equation}
  Then $T(L)=K$. By \cref{lem:linear-operator}, $\Gamma_l(K)=\Gamma_l(T(L))\leq \Gamma_l(K_m^{1*})$.

  (b). Let $\pm v_1,\dots,\pm v_m$ be the vertices of $K$. Put
  $u_i=(v_i,e_i^m),~\forall i\in[m]$. Clearly, $u_1,\dots,u_m$ are linearly
  independent. Then $L:=\CH\qty{\pm u_1,\dots,\pm u_m}$ is an $m$-dimensional
  crosspolytope and $K=T(L)$, where $T$ is defined by \cref{eq:map}. From
  \cref{lem:linear-operator}, it follows that $\Gamma_l(K)\leq
  \Gamma_l(L)=\Gamma_l(K_m^1)$.
\end{proof}

\begin{thm}
  \label{thm:simplex}
  Let $K\subseteq \R^n$ be a convex polytope.
  \begin{enumerate}
  \item[(a).] If $K$ has $m$ $(m\geq n+1)$ vertices and $t=2^{\frac{n}{m-1}}$, then
    \begin{displaymath}
      \Gamma_{2^n}(K)\leq\frac{m-1}{m-1+k(m-1,t)}\leq \frac{m-1}{m-1+\floor{a\qty(t)(m-1)}}.
    \end{displaymath}
  \item[(b).] If $K$ is centrally symmetric and has $2m$ $(m\geq n)$ vertices,
    and $t=2^{\frac{n}{m}}$, then
    \begin{displaymath}
      \Gamma_{2^n}(K)\leq \frac{m}{m+l(m,t)}\leq \frac{m}{m+\floor{b\qty(t)m}}.
    \end{displaymath}
  \end{enumerate}  
\end{thm}
\begin{proof}
  (a). Since $m>n>0$, we have $t\in (1,2]$. By \cref{thm:simplex-crosspolytope} and \cref{prop:Kpn}, we have
  \begin{align*}
    \Gamma_{2^n}(K)\leq\Gamma_{2^n}(K_{m-1}^{1*})=\Gamma_{t^{m-1}}(K_{m-1}^{1*})&\leq \frac{m-1}{m-1+k(m-1,t)}\\
                                                             &\leq \frac{m-1}{m-1+\floor{a(t)(m-1)}}.
  \end{align*}

  (b). Since $m\geq n$, $1<t\leq 2$. By \cref{thm:simplex-crosspolytope} and \cref{prop:Kpn}, we have
  \begin{displaymath}
    \Gamma_{2^n}(K)\leq\Gamma_{2^n}(K_m^1)\leq \frac{m}{m+l(m,t)}\leq\frac{m}{m+\floor{b\qty(2^{\frac{n}{m}})m}}.\qedhere
  \end{displaymath}
\end{proof}
\begin{rem}
  \cref{thm:simplex} and (b) in \cref{lem:estimations-k(n)} show that, when
  $m\leq t^{m-1}=2^n$, we have
  \begin{displaymath}
    \Gamma_{2^n}(K)\leq \frac{m-1}{m}<1.
  \end{displaymath}
  In particular, when $m=2^n$,
  \begin{displaymath}
    \Gamma_{2^n}(K)\leq \frac{2^n-1}{2^n-1+1}=1-\frac{1}{2^n},
  \end{displaymath}
  which is worse than \cref{eq:general-estimation}. When $m-1$ is proportional
  to $n$, (a) in \cref{thm:simplex} yields better estimations when $n$ is large.
  Take the case when $m-1=2n$ for example. We have
  \begin{displaymath}
    \Gamma_{2^n}(K)\leq \frac{m-1}{m-1+\floor{a(\sqrt{2})(m-1)}}\approx
    \frac{2n}{2n+\floor{0.104828 \cdot(2n)}}.
  \end{displaymath}
  It follows that
  \begin{displaymath}
    \limsup\limits_{n\to\infty}\Gamma_{2^n}(K)\leq \frac{1}{1.104828}\approx 0.905118,
  \end{displaymath}
  which is much better than the general estimation \cref{eq:general-estimation} for large $n$.
  \cref{thm:simplex} (a) is also valid when $n$ is small. Take the case when
  $n=6$ and $m=9$ for example. In this case $t=2^{\frac{6}{8}}$ and
  $k(m-1,t)=2$. It follows that
  \begin{displaymath}
    \Gamma_{64}(K)\leq \frac{8}{8+2}=0.8
  \end{displaymath}
  holds for each $K\in\SCB[6]$ having $9$ vertices. This is also better than
  \cref{eq:general-estimation}, which gives $\Gamma_{64}(K)\leq 1-\frac{1}{7}\approx 0.86$.

  For (b) in \cref{thm:simplex}, the situation is similar.
\end{rem}

\section{Acknowledgement}
The authors are supported by the National Natural Science Foundation of China
(grant numbers 12071444 and 12001500), the Natural Science Foundation of
Shanxi Province of China (grant number 201901D111141), and the Scientific and
Technological Innovation Programs of Higher Education Institutions in Shanxi
(grant numbers 2020L0290 and 2020L0291).



\begin{thebibliography}{10}

\bibitem{Betke-Henk1993}
U.~Betke and M.~Henk, \emph{Intrinsic volumes and lattice points of
  crosspolytopes}, Monatsh. Math. \textbf{115} (1993), no.~1-2, 27--33.
  \MR{1223242}

\bibitem{Bezdek2010}
K.~Bezdek, \emph{Classical {T}opics in {D}iscrete {G}eometry}, CMS Books in
  Mathematics/Ouvrages de Math\'ematiques de la SMC, Springer, New York, 2010.
  \MR{2664371 (2011j:52014)}

\bibitem{Bezdek-Khan2018}
K.~Bezdek and Muhammad~A. Khan, \emph{The geometry of homothetic covering and
  illumination}, Discrete Geometry and Symmetry, Springer Proc. Math. Stat.,
  vol. 234, Springer, Cham, 2018, pp.~1--30. \MR{3816868}

\bibitem{Boltyanski-Martini-Soltan1997}
V.~Boltyanski, H.~Martini, and P.S. Soltan, \emph{Excursions into
  {C}ombinatorial {G}eometry}, Universitext, Springer-Verlag, Berlin, 1997.
  \MR{1439963}

\bibitem{Brass-Moser-Pach2005}
P.~Brass, W.~Moser, and J.~Pach, \emph{Research {P}roblems in {D}iscrete
  {G}eometry}, Springer, New York, 2005. \MR{2163782 (2006i:52001)}

\bibitem{Martini-Soltan1999}
H.~Martini and V.~Soltan, \emph{Combinatorial problems on the illumination of
  convex bodies}, Aequationes Math. \textbf{57} (1999), no.~2-3, 121--152.
  \MR{1689190 (2000b:52006)}

\bibitem{Polya-Szego1998}
G.~P\'{o}lya and G.~Szeg\H{o}, \emph{Problems and theorems in analysis. {I}},
  Classics in Mathematics, Springer-Verlag, Berlin, 1998, Series, integral
  calculus, theory of functions, Translated from the German by Dorothee Aeppli,
  Reprint of the 1978 English translation. \MR{1492447}

\bibitem{Robbins1955}
H.~Robbins, \emph{A remark on {S}tirling's formula}, Amer. Math. Monthly
  \textbf{62} (1955), 26--29. \MR{69328}

\bibitem{Jaegermann1989}
N.~Tomczak-Jaegermann, \emph{Banach-{M}azur {D}istances and
  {F}inite-dimensional {O}perator {I}deals}, Pitman Monographs and Surveys in
  Pure and Applied Mathematics, vol.~38, Longman Scientific \& Technical,
  Harlow; copublished in the United States with John Wiley \& Sons, Inc., New
  York, 1989. \MR{993774}

\bibitem{Wu-He2019}
Senlin Wu and Chan He, \emph{Covering functionals of convex polytopes}, Linear
  Algebra Appl. \textbf{577} (2019), 53--68. \MR{3943803}

\bibitem{Xue-Lian-Zhang2021+}
Fei Xue, Yanlu Lian, and Yuqin Zhang, \emph{On {H}adwiger's covering functional
  for the simplex and the cross-polytope}, 2021.

\bibitem{Zong2010}
Chuanming Zong, \emph{A quantitative program for {H}adwiger's covering
  conjecture}, Sci. China Math. \textbf{53} (2010), no.~9, 2551--2560.
  \MR{2718847 (2012c:52040)}

\end{thebibliography}

\def\polhk#1{\setbox0=\hbox{#1}{\ooalign{\hidewidth
  \lower1.5ex\hbox{`}\hidewidth\crcr\unhbox0}}} \def\cprime{$'$}
  \def\polhk#1{\setbox0=\hbox{#1}{\ooalign{\hidewidth
  \lower1.5ex\hbox{`}\hidewidth\crcr\unhbox0}}} \def\cprime{$'$}
\providecommand{\bysame}{\leavevmode\hbox to3em{\hrulefill}\thinspace}
\providecommand{\MR}{\relax\ifhmode\unskip\space\fi MR }
\providecommand{\MRhref}[2]{%
  \href{http://www.ams.org/mathscinet-getitem?mr=#1}{#2}
}
\providecommand{\href}[2]{#2}

\end{document}